\documentclass[12pt, a4paper, reqno]{amsart}
\usepackage[left=2.5cm,right=2.2cm,top=2.5cm,bottom=3.5cm]{geometry}
\usepackage{graphicx}
\usepackage[T2A]{fontenc}
\usepackage[utf8]{inputenc}
\usepackage[english]{babel}
\usepackage{amsmath, amssymb, amsthm, mathtools}
\usepackage[hidelinks]{hyperref}
\usepackage{tikz-cd, quiver}
\usetikzlibrary{babel}

\title{Unramified case of Colliot-Th{\'e}l{\`e}ne conjecture on quadratic forms}
\author{Ivan Panin, Ivan Gaidai-Turlov}
\date{}

\makeatletter 
\thm@space@setup
\newcommand{\ps}{\par\addvspace{\thm@preskip}}
\makeatother

\newcommand{\ul}{\underline}
\newcommand{\AffSp}{\mathbb{A}}           
\newcommand{\ProjSp}{\mathbb{P}}          

\newcommand{\red}{\mathrm{red}}           
\newcommand{\rat}{\mathcal R}             

\newcommand{\basedvr}{R}

\DeclareMathOperator{\spec}{Spec}

\DeclareMathOperator{\nil}{Nil}

\newtheorem{theorem}{Theorem}[section]
\newtheorem{lemma}[theorem]{Lemma}
\newtheorem{corollary}[theorem]{Corollary}
\newtheorem{conjecture}[theorem]{Conjecture}
\theoremstyle{definition}
\newtheorem{construction}[theorem]{Construction}
\newtheorem{remark}[theorem]{Remark}
\newtheorem{definition}[theorem]{Definition}

\addto\captionsenglish{}

\counterwithin{theorem}{subsection}
\renewcommand{\thetheorem}{%
  \ifnum\value{subsection}>0
    \thesubsection.%
  \else
    \thesection.%
  \fi
  \arabic{theorem}%
}

\newcommand{\widebar}[1]{\overline{#1}} 
\makeatletter
\newcommand*\rel@kern[1]{\kern#1\dimexpr\macc@kerna}
\renewcommand*\widebar[1]{%
  \begingroup
  \def\mathaccent##1##2{%
    \rel@kern{0.8}%
    \smash{\overline{\rel@kern{-0.8}\macc@nucleus\rel@kern{0.2}}}%
    \rel@kern{-0.2}%
  }%
  \macc@depth\@ne
  \let\math@bgroup\@empty \let\math@egroup\macc@set@skewchar
  \mathsurround\z@ \frozen@everymath{\mathgroup\macc@group\relax}%
  \macc@set@skewchar\relax
  \let\mathaccentV\macc@nested@a
  \macc@nested@a\relax111{#1}%
  \endgroup
}
\makeatother

\begin{document}

\begin{abstract}
    The Colliot-Th{\'e}l{\`e}ne conjecture on quadratic forms asserts that any rationally isotropic nondegenerate quadratic form over a semilocal regular domain with 2 invertible is isotropic. We prove the Colliot-Th{\'e}l{\`e}ne conjecture for semilocal regular domains that are geometrically regular over a discrete valuation ring with 2 invertible. In particular, this includes unramified regular local rings of mixed characteristic $(0,p)$, where $p\neq 2$.
\end{abstract}

\footnote{This work was supported by the Ministry of Science and Higher Education of the Russian Federation (agreement 075-15-2025-344 dated 29/04/2025 for Saint Petersburg Leonhard Euler International Mathematical Institute at PDMI RAS).}

\maketitle

\section{Introduction}

This paper further extends the results of \cite{PaninPimenov2016, PaninTyurin2024}. See also \cite{PaninPimenov2010, Panin2009}. Our main goal is to prove the Colliot-Th{\'e}l{\`e}ne conjecture on quadratic forms \cite{ColliotThelene1979} for regular semilocal domains geometrically regular over a dvr, covering in particular the unramified mixed characteristic case. Our main result is Theorem \ref{ct-smooth} below.

\begin{conjecture}[Colliot-Th{\'e}l{\`e}ne]\label{ct-regular}
    Let $\mathcal O$ be a semilocal regular domain with $2 \in \mathcal O^\times$. Any rationally isotropic (i.e., isotropic over the fraction field $K = \rat(\mathcal O)$) nondegenerate quadratic form $q$ over $\mathcal O$ is isotropic.
\end{conjecture}

\begin{theorem}\label{ct-smooth}
    Let $\mathcal O$ be a semilocal essentially smooth domain over a dvr $\basedvr$ with $2 \in \basedvr^\times$. Any rationally isotropic nondegenerate quadratic form $q$ over $\mathcal O$ is isotropic.
\end{theorem}

The following corollaries are consequences of our main theorem.

\begin{corollary}\label{ct-geom-regular}
    Let $\mathcal O$ be a semilocal regular domain geometrically regular over a dvr $\basedvr$ with $2 \in \basedvr^\times$. Any rationally isotropic nondegenerate quadratic form $q$ over $\mathcal O$ is isotropic.
\end{corollary}
\begin{proof}
    The corollary is deduced from Theorem \ref{ct-smooth} by the same argument as in the proof of \cite[Theorem 1.1]{PaninPimenov2010}. The only difference is that our base is a dvr rather than a field. We use Popescu's theorem \cite[Theorem 2.5]{Popescu1986}, which requires the additional assumption of geometric regularity. See also Swan's exposition \cite{Swan1998}.
\end{proof}

\begin{corollary}
    Let $\mathcal O$ be a semilocal regular domain geometrically regular over a dvr $\basedvr$ with $2\in\basedvr^\times$. Let $q$ be a nondegenerate quadratic form over $\mathcal O$, and let $a\in\mathcal O^\times$. If the equation $q=a$ has a solution over $K = \rat(\mathcal O)$, then it has a solution over $\mathcal O$.
\end{corollary}
\begin{proof}
    This follows from Corollary \ref{ct-geom-regular} and \cite[Proposition~1.2(b)]{ColliotThelene1979}.
\end{proof}

To state the next corollary, we recall the notion of an unramified quadratic form. A nondegenerate quadratic form $q_K$ over the fraction field $K = \rat(\mathcal O)$ of a local integral domain $\mathcal O$ is called unramified if, for every height-one prime ideal $\mathfrak p$ of $\mathcal O$, there exists a nondegenerate quadratic form $q_{\mathfrak p}$ over the localization $\mathcal O_{\mathfrak p}$ such that $q_{\mathfrak p} \otimes_{\mathcal O_{\mathfrak p}} K \cong q_K$.

\begin{corollary}[Purity]
    Let $\mathcal O$ be a regular local domain geometrically regular over a dvr $\basedvr$ with $2\in\basedvr^\times$. Let $q_K$ be a nondegenerate quadratic form over $K = \rat(\mathcal O)$ that is unramified over $\mathcal O$. Then there exists a nondegenerate quadratic form $q$ over $\mathcal O$ such that $q \otimes_{\mathcal O} K \cong q_K$.
\end{corollary}
\begin{proof}
    The corollary follows from Corollary \ref{ct-geom-regular} by the same argument as in the proof of \cite[Corollary 3.1]{PaninPimenov2010}, using the purity theorem for Witt groups over regular local rings geometrically regular over a dvr \cite[Theorem 8.2 and the subsequent corollary]{GillePanin2024}.
\end{proof}

\begin{remark}
    Recall that a regular local ring $\mathcal O$ of mixed characteristic $(0,p)$ is called unramified if $\mathcal O/p\mathcal O$ is regular. Equivalently, $\mathcal O$ is geometrically regular over $\mathbb Z_{(p)}$. In particular, if $p\neq 2$, such rings satisfy the assumptions of Corollary \ref{ct-geom-regular}.
\end{remark}

It is worth mentioning that the Colliot-Th{\'e}l{\`e}ne conjecture is a special case of a more general conjecture on projective homogeneous schemes, sometimes attributed to Colliot-Th{\'e}l{\`e}ne or Panin. In the \cite{PaninPimenov2010} one can find historical remarks. The conjecture is stated in \cite{Cesnavicius2022} as follows.

\begin{conjecture}
    Let $\mathcal O$ be a regular local ring and let $G$ be a reductive $\mathcal O$-group scheme. If the generic fiber of $G$ has a parabolic subgroup, then $G$ itself has a parabolic subgroup of the same type.
\end{conjecture}

The conjecture remains open in general, even for regular local rings essentially smooth over a field. In addition to the results on the Colliot-Th{\'e}l{\`e}ne conjecture mentioned above, the case of Borel subgroups over regular local rings containing a field is settled in \cite[Theorem 9.5]{Cesnavicius2022}, the case of inner forms of type $A_n$ in characteristic zero follows from \cite{Ojanguren2015}, and certain cases involving exceptional groups are treated in \cite{PaninPetrov2015}.

\ps

The proof of our main theorem, Theorem \ref{ct-smooth}, relies on three results: geometric Lemma \ref{geom-lemma}, which is a suitable version of \cite[Lemma 10.1]{OjangurenPanin2001}, preparatory Lemma \ref{generic-nondegeneracy}, and the following version of Springer's theorem due to Gille and Neher.

\begin{theorem}[{\cite[Theorem 2.1]{GilleNeher2021}}]\label{springer}
    Let $A$ be a semilocal ring with $2 \in A^\times$ and let $q$ be a nondegenerate quadratic form over $A$. Let $B$ be a finite projective one-generated $A$-algebra of constant odd degree. If $q$ is $B$-isotropic, then $q$ is $A$-isotropic.
\end{theorem}

We rely heavily on two features of Theorem \ref{springer}: it does not require the residue fields to be infinite, and it assumes only that the algebra $B$ is one-generated, rather than etale. The latter feature is used both in the reduction to perfect residue fields and in geometric Lemma \ref{geom-lemma}, where we only ensure that $D_0$ and $D_1$ are one-generated over the base, since etaleness cannot be expected in general.

\ps

Throughout the paper, we fix the following notation. Let $\basedvr$ be a discrete valuation ring of mixed characteristic $(0,p)$, $p\neq 2$, and let $V=\spec\basedvr=\{v,\eta\}$, where $v$ is the closed point and $\eta$ is the generic point. For any $V$-schemes $X$ and $S$, we denote by $X_S$ the $S$-scheme obtained from $X$ by base change along $S\to V$. In particular, $X_v$ and $X_\eta$ denote the closed and generic fibers of $X$. By a closed or generic fiber we mean a fiber over a closed or generic point of the base, respectively. For an integral scheme $X$, we denote by $\rat(X)$ its field of rational functions. By an essentially smooth algebra over a ring $A$, we mean an $A$-algebra that is a localization of a smooth $A$-algebra. For a $V$-scheme $X$, we denote its regular and smooth locus by $X_{\operatorname{reg}}$ and $X_{\operatorname{sm}}$, respectively. If $X$ is locally of finite type over $V$, then $X_{\operatorname{reg}}$ is open, since dvrs of mixed characteristic are excellent. A scheme is called semilocal if it is the spectrum of a semilocal ring. We say that an $S$-scheme $X$ is of constant relative dimension $n$ if every fiber $X_s$, $s\in S$, is nonempty and has dimension $n$.

\section{Reduction to the case of an infinite perfect residue field}

\begin{lemma}\label{dvr-extension}
    Let $R$ be a dvr with uniformizer $\pi$ and residue field $k$. Let $k'$ be a finite simple extension of $k$. Then there exists a dvr $R'$, finite, flat and one-generated over $R$, with residue field $k'$ and uniformizer $\pi$.
\end{lemma}
\begin{proof}
    Write $k' = k[t]/(f)$ for some monic irreducible polynomial $f \in k[t]$. Choose a monic lift $F \in R[t]$ of $f$, and set $R' = R[t]/(F)$. Since $F$ is monic, $R'$ is finite flat over $R$. Moreover, $R'/\pi R' \cong k'$, so $\pi R'$ is the unique maximal ideal of $R'$. Thus $R'$ is a one-dimensional regular local ring and hence a dvr.
\end{proof}

\begin{construction}\label{dvr-tower}
    Let $R$ be a dvr with uniformizer $\pi$ and residue field $k$. Denote $p = \operatorname{char}k$.
    \begin{enumerate}
        \item Assume that $k$ is finite. Fix an algebraic closure of $k$. For each $\alpha\in\mathbb Z_{\geqslant 0}$, let $k^\alpha$ be the unique extension of $k$ of degree $p^{2^\alpha}$ inside this closure. Starting with $R^0=R$, we recursively apply Lemma~\ref{dvr-extension} to construct dvrs $R^\alpha$ with residue fields $k^\alpha$ and uniformizer $\pi$ such that each $R^\alpha$ is finite, flat and one-generated over $R^{\alpha-1}$.

        \item Assume that $k$ is imperfect, and denote its perfect closure by $k^{\mathrm{perf}}$. Choose a well-ordering $(a_\alpha)_{\alpha<\lambda}$ of $k^{\mathrm{perf}}$. For each ordinal $\alpha\leq\lambda$, let $k^\alpha\subset k^{\mathrm{perf}}$ be the subfield generated by $k$ and all elements $a_\beta$ with $\beta<\alpha$. We construct dvrs $R^\alpha$ with residue fields $k^\alpha$ and uniformizer $\pi$ by transfinite induction on $\alpha$. Set $R^0=R$. At a successor stage $\alpha=\beta+1$, the extension $k^\alpha/k^\beta$ is finite and simple of $p$-primary degree. Lemma~\ref{dvr-extension} then yields a dvr $R^\alpha$ with residue field $k^\alpha$ and uniformizer $\pi$ that is finite, flat and one-generated over $R^\beta$. At a limit ordinal $\alpha$, define $R^\alpha$ as the direct limit $\bigcup_{\beta<\alpha}R^\beta$. Since all transition maps preserve the uniformizer $\pi$, the ring $R^\alpha$ is again a dvr with uniformizer $\pi$.
    \end{enumerate}
    Let $R^\infty$ be the direct limit $\bigcup_{\alpha} R^\alpha$ of the constructed tower, and denote its residue field by $k^\infty$. By construction, $R^\infty$ is a dvr with uniformizer $\pi$, and $k^\infty$ is infinite and perfect.
\end{construction}

Using this construction and a version of Springer's Theorem~\ref{springer}, we can reduce Theorem \ref{ct-smooth} to the case of a dvr with infinite perfect residue field.

\begin{lemma}\label{infinite-perfect-reduction}
Assume that Theorem \ref{ct-smooth} holds for every dvr $R$ with infinite perfect residue field and $2$ invertible. Then it holds for every dvr with $2$ invertible.
\end{lemma}
\begin{proof}
    Let $\mathcal O$ and $q$ be as in Theorem \ref{ct-smooth}. Assume that the residue field $k$ of $R$ is finite or imperfect. Consider $\mathcal O^\alpha=\mathcal O\otimes_R R^\alpha$ and $\mathcal O^\infty=\mathcal O\otimes_R R^\infty$, where $R^\infty$ is as in Construction~\ref{dvr-tower}. The form $q$ remains nondegenerate and rationally isotropic over $\mathcal O^\infty$. Hence, by assumption, it is $\mathcal O^\infty$-isotropic.

    We can then descend using Theorem~\ref{springer}. At a successor stage $\alpha=\beta+1$, descent is possible since the extension $\mathcal O^\alpha/\mathcal O^\beta$ is finite, flat and one-generated of constant odd degree, as $\operatorname{char}k\neq 2$. At a limit ordinal $\alpha$, an isotropic unimodular vector is already defined over some $\mathcal O^\beta$ with $\beta<\alpha$, where it is also isotropic and unimodular by faithfully flat descent.
\end{proof}

\section{Auxiliary results}

\begin{lemma}\label{one-generated-algebras-k}
    A finite algebra over an infinite perfect field is one-generated if and only if its nilradical is (locally) principal.
\end{lemma}
\begin{proof}
    Obviously, the nilradical of a finite one-generated algebra is locally principal. Let us prove the converse. Denote the field by $k$ and the algebra by $A$. Since $A$ is artinian, we can decompose $A$ into the product of local algebras $A_i$. 
    
    Assume that each algebra $A_i$ is generated by an element $a_i \in A_i$. Since $k$ is infinite, we can add to each generator some element $b_i \in k$ such that new generators $a_i + b_i$ have pairwise coprime minimal polynomials. By the Chinese remainder theorem, element $(a_i+b_i)_i$ generates the product $A = \prod A_i$.

    It remains to consider the case of a local algebra $A$. In this case, the nilradical of $A$ is maximal and generated by an element $\varepsilon$. Since $k$ is perfect, the field $A/\varepsilon A$ is a finite separable extension of $k$, and hence is generated by an element $\bar a \in A/\varepsilon A$. Choose a lift $a \in A$ of $\bar a$, and let $f$ be the minimal polynomial of $\bar a$, which is separable. Modulo $\varepsilon^2 A$, we have $f(a+\varepsilon b) \equiv f(a)+\varepsilon b f'(a)$, and $f'(a)$ is invertible. For a suitable choice of $b \in A$, the element $f(a+\varepsilon b)$ generates the nilpotent ideal $\varepsilon A$. Since $a+\varepsilon b$ also generates $A/\varepsilon A$, we conclude that $a+\varepsilon b$ generates $A$.
\end{proof}

\begin{remark}
    Such algebras are precisely the finite products of algebras of the form $L[\varepsilon]/(\varepsilon^d)$, where $L/k$ is a finite field extension.
\end{remark}

\begin{lemma}\label{one-generated-algebras}
    Let $S$ be a semilocal scheme, and let $D$ be a finite $S$-scheme. Assume that the residue fields $\rat(s)$ of the closed points $s \in S$ are infinite and perfect. Suppose that $\nil(D_s)$ is locally principal for every closed point $s \in S$. Then $D$ is one-generated over $S$.
\end{lemma}
\begin{proof}
    By Lemma \ref{one-generated-algebras-k}, each $D_s$ is one-generated over $s$. By the Chinese remainder theorem, we can choose a common lift of these
    generators. This lift generates $D$ over $S$ by Nakayama's lemma.
\end{proof}

\begin{lemma}\label{fiber-product-nilradical}
    Let $X$ be a scheme, and let $Y$ and $Z$ be $X$-schemes. Assume that the fiber product $Z \times_X Y_\red$ is reduced. Then $\nil(\mathcal O_{Z \times_X Y}) = \nil(\mathcal O_Y)\mathcal O_{Z \times_X Y}$.

    In particular, if $\nil(\mathcal O_Y)$ is locally principal on $Y$, then $\nil(\mathcal O_{Z \times_X Y})$ is locally principal on $Z \times_X Y$ as well.
\end{lemma}
\begin{proof}
    The ideal sheaf $\nil(\mathcal O_Y)\mathcal O_{Z \times_X Y}$ is a nil ideal and defines the closed reduced subscheme $Z \times_X Y_\red$ in $Z \times_X Y$. Therefore, it coincides with $\nil(\mathcal O_{Z \times_X Y})$.
\end{proof}

\begin{definition}\label{adm-comp}
    Let $S$ be a scheme. We say that an $S$-scheme $X$ admits a compactification if there exists a projective $S$-scheme $\widebar X$ such that $X$ is a fiberwise dense open subscheme of $\widebar X$ and the complement $X^\infty = \widebar X - X$ is the zero locus of a section $f_0 \in \Gamma(\widebar X, \mathcal O_{\widebar X}(1))$ of a very ample over $S$ invertible sheaf $\mathcal O_{\widebar X}(1)$ on $\widebar X$. In such a case we say that $\widebar X$ is a compactification of $X$ over $S$.
\end{definition}

\begin{lemma}\label{adm-comp-vs-normal}
    Let $S$ be a semilocal normal scheme and let $X$ be an irreducible $S$-scheme. The following are equivalent
    \begin{enumerate}
        \item $X$ admits a compactification over $S$ and $\dim X_s=n$ for every $s\in S$.
        \item There exists a finite surjective $S$-morphism $X \to \AffSp^n_S$
    \end{enumerate}
\end{lemma}
\begin{proof}
    Assume (1). A finite surjective morphism can be constructed by choosing suitable sections of the very ample invertible sheaf $\mathcal O_{\widebar X}(d)$ for sufficiently large $d$. 
    
    Assume (2). By Zariski's Main Theorem applied to the composition $X \to \AffSp^n_S \to \ProjSp^n_S$, there exists a finite morphism $\pi\colon \widebar X \to \ProjSp^n_S$ such that $X \subset \widebar X$ is an open subscheme. Replacing $\widebar X$ by the scheme-theoretic closure of $X$ in $\widebar X$, we can achieve that $X = \pi^{-1}(\AffSp^n_S)$. Since $\pi$ is finite and dominant, it is surjective. Moreover, since $\ProjSp^n_S$ is normal, the finite surjective morphism $\widebar X_{\mathrm{red}}\to\ProjSp^n_S$ satisfies going down, implying fiberwise density of $X$ in $\widebar X$.
\end{proof}

\begin{remark}
    In Lemma \ref{adm-comp-vs-normal}, the dimension condition in (1) can be strengthened to require that every fiber $X_s$ be nonempty and of pure dimension $n$. Indeed, under (2), going down implies that every irreducible component of $X_s$ dominates $\AffSp^n_{s}$ and hence has dimension $n$.
\end{remark}

\begin{remark}
    If, in addition, $X$ is Cohen--Macaulay and $S$ is regular, then the finite surjective morphism in (2) is flat by miracle flatness. In particular, $X$ is flat over $S$.
\end{remark}

\section{Geometric lemma}

\begin{theorem}\label{geom-lemma}
    Let $S$ be a semilocal regular scheme. Let $X$ be an irreducible regular $S$-scheme of constant relative dimension $n \geqslant 1$ that admits a compactification. Let $\Delta\colon S \to X$ be a section such that $X$ is smooth over $S$ along $\Delta(S)$, and let $Z \subset X$ be a closed subset such that the complement $X-Z$ is dense in each closed fiber $X_s$ of $X$ over $S$. Assume that
    \begin{itemize}
        \item The residue fields $\rat(s)$ of the closed points $s \in S$ are infinite and perfect.
        \item The nilradical $\nil(X_s)$ of each closed fiber is generically principal.
    \end{itemize}
    Then there exists a finite flat $S$-morphism
    \[
    \pi\colon X \to \AffSp^n_S
    \]
    such that
    \begin{enumerate}
        \item $\pi^{-1}(0_S) = \Delta(S) \sqcup D_0$ and $D_0 \cap Z = \varnothing$;
        \item $\pi^{-1}(1_S) = D_1$ and $D_1 \cap Z = \varnothing$;
        \item $D_1$ and $D_0$ are finite, flat, and one-generated over $S$.
    \end{enumerate}
    Here $0_S = (0, 0, \dots, 0) \in \AffSp_S^n(S)$ and $1_S = (1, 1, \dots, 1) \in \AffSp_S^n(S)$ are $S$-points.
\end{theorem}
\begin{proof}
    Let $\widebar X$, $\mathcal O_{\widebar X}(1)$ and $f_0$ be as in Defenition \ref{adm-comp}. Let $\widebar Z$ denote the closure of $Z$ in $\widebar X$. 
    
    First, we solve the problem on the closed fiber $X_s$ for each closed point $s \in S$. Note that such a fiber may be nonreduced. Thus, we consider the corresponding reduced scheme $X_{s, \red}$. Since $\rat(s)$ is perfect, the smooth locus of $\widebar X_{s, \red}$ over $s$ is open and dense. Moreover, by assumption, $\nil(X_s)$ is locally principal on some dense open subset. Note that $\nil(X_s)$ vanishes at $\Delta(s)$ since $X_s$ is smooth over $s$ at $\Delta(s)$. We choose a closed subset $Y_s \subset \widebar X_s$ such that $\Delta(s) \notin Y_s$, $\widebar X_{s, \red} - Y_s$ is smooth over $s$, and $\nil(\widebar X_s - Y_s)$ is locally principal. Note that $\dim Y_s \leqslant n-1$ and $\dim \widebar Z_s \leqslant n-1$.

    We construct a finite surjective morphism $\widebar\pi \colon \widebar X \to \ProjSp^n_S$ by choosing sections $f_i$ of the sheaf $\mathcal{O}_{\widebar{X}}(2)$, $i = 1, \dots, n$:
    \[
    \widebar\pi = \left[f_0^2 : f_1 : f_2 : \dots : f_n\right] \colon \widebar X \to \ProjSp^n_S.
    \]
    The desired morphism $\pi$ will be the restriction of $\widebar\pi$ to the open subset $X = \{f_0 \neq 0\}$.

    We can choose such sections on each closed fiber $\widebar X_s$, one by one, using a version of Bertini's theorem. Since the field $\rat(s)$ is infinite, at step $r$ we can choose a section $f_r$ of $\mathcal{O}_{\widebar{X}_s}(2)$ such that the subscheme $\{f_r = 0\}$
    \renewcommand{\theenumi}{\roman{enumi}}%
    \begin{enumerate}
        \item contains $\Delta(s)$,
        \item intersects $\bigcap_{i = 1}^{r-1} \{f_i = 0\}$ transversely in the smooth locus of $\widebar X_{s, \red}$ over $s$,
        \item avoids a point in each irreducible component of $\bigcap_{i = 1}^{r-1} \{f_i = 0\} \bigcap X^\infty_s$,
        \item avoids a point in each irreducible component of $\bigcap_{i = 1}^{r-1} \{f_i = 0\} \bigcap Y_s$,
        \item avoids a point in each irreducible component of $(\bigcap_{i = 1}^{r-1} \{f_i = 0\} \bigcap \widebar Z_s) - \Delta(s)$.
    \end{enumerate}
    We then choose simultaneous lifts of these sections such that $\Delta(S) \subset \{f_i = 0\}$. This is possible by property $(\mathrm{i})$.

    Property $(\mathrm{iii})$ ensures that the morphism $\widebar\pi$ is everywhere defined. Indeed, the common zero locus of $f_0^2, f_1, \dots, f_n$ is a closed subset with empty closed fibers. Since $\widebar X$ is projective over $S$, this locus is empty. Note that $\widebar\pi$ is affine and projective, hence finite. Since both $X$ and $\AffSp^n_S$ are regular, the morphism $\pi$ is flat by miracle flatness.  By dimension considerations, $\widebar\pi$ is surjective.

    By construction, $\Delta$ is a section of the finite flat morphism $\pi^{-1}(0_S) \to 0_S$. We claim that this morphism is etale along $\Delta(S)$, implying that $\Delta(S)$ is an open and closed subscheme of $\pi^{-1}(0_S)$. It suffices to check etaleness after base change to each closed point $s \in S$. 
    Property $(\mathrm{i})$ ensures that $\pi^{-1}(0_s)$ meets the closed fiber $\widebar X_{s}$ at $\Delta(s) \in \widebar X_{s, \mathrm{sm}}$. By property $(\mathrm{ii})$, this intersection is transverse. Hence $\pi^{-1}(0_s)$ is smooth over $s$ at $\Delta(s)$, as required. We write $\pi^{-1}(0_S) = \Delta(S) \sqcup D_0$, as in the statement of the theorem.

    Property $(\mathrm{v})$ ensures that $D_ 0\cap Z$ is a closed subset of $D_0$ with empty closed fibers and is therefore empty, i.e., $D_ 0\cap Z = \varnothing$.

    By properties $(\mathrm{ii})$ and $(\mathrm{iv})$, the intersection $D_0 \cap \widebar X_{s,\red}$ is smooth over $s$ and, in particular, reduced. Hence, by Lemma \ref{fiber-product-nilradical}, the nilradical $\nil(D_{0,s})$ is locally principal. Therefore, by Lemma \ref{one-generated-algebras}, $D_0$ is one-generated over $S$.

    To ensure the required properties of $D_1$, choose a rational point $p_s$ in each closed fiber $\AffSp^n_s$ that lies in the open subset over which $\pi|_{X_{s, \red}}\colon X_{s, \red} \to \AffSp^n_s$ is etale and avoids $0_s$, $\widebar\pi(\widebar Z_s)$, and $\widebar\pi(Y_s)$. Note that, by construction, $\pi|_{X_{s,\red}}$ is etale over $0_s$, so it is etale over some open neighborhood of $0_s$. Since the residue field $\rat(s)$ is infinite, such a point exists. We then choose a common lift $p$ of these points and make a linear change of coordinates in $\AffSp^n_S$ sending $p$ to $1_S$. By the same argument as for $D_0$, the subscheme $D_1$ satisfies all the required properties.
\end{proof}

\begin{remark}
    We use a slightly modified version of Bertini's theorem \cite[Theorem 2.1]{Artin1973Bertini}. The additional requirement of avoiding finitely many prescribed points can be imposed since, for any two distinct points, the set of sections vanishing at the first point but not at the second is open and dense in the space of sections vanishing at the first point.
\end{remark}

\begin{remark}
    The regularity assumption on $X$ can be weakened to Cohen--Macaulayness at no additional cost. This is useful when constructing schemes $X$ admitting a compactification over non-excellent bases, such as certain discrete valuation rings of positive characteristic.
\end{remark}

\section{Proof of Theorem \ref{ct-smooth}}

By Lemma \ref{infinite-perfect-reduction} we can and will assume that the residue field of $R$ is infinite and perfect. The case of the Colliot-Th{\'e}l{\`e}ne conjecture for $\dim(\mathcal O) = 1$ is well-known, so we will assume that $\dim(\mathcal O) \geqslant 2$. By assumption, $\mathcal O$ is essentially smooth $\basedvr$-algebra. Hence, we can choose an integral normal projective $V$-scheme $\widebar X$ of constant relative dimension and a finite set of points $\ul x \subset \widebar X_{\operatorname{sm}}$ such that $\mathcal O \cong \mathcal O_{\ul x}$. For example, choose a smooth affine $V$-scheme $X_0$ with a quasi-finite dominant morphism $X_0 \to \AffSp^n_V$ and normalize $\ProjSp^n_V$ in the field of fractions $\rat(X_0)$. Note that the normalization is finite, since dvrs of mixed characteristic are excellent. Set $U = \spec \mathcal O_{\ul x}$.

Let $\ul\theta$ be the set of generic points $\theta_i \in \widebar X_v$. The ring $\mathcal O_{\theta_i}$ is a dvr since it is normal and one-dimensional. We denote by $\pi_i$ a uniformizer of $\mathcal O_{\theta_i}$ that is invertible in $\mathcal O_{\theta_j}$ for $j \neq i$, and by $v_i$ the corresponding valuation.

We can diagonalize our form and write $q = \sum a_kT_k^2$, since $\mathcal O_{\ul x}$ is semilocal with $2$ inverted. Note that $v_i(a_k) = 0$ if $x \in \widebar{\{\theta_i\}}$ for some $x \in \ul x$, since the form $q$ is nondegenerate over $\mathcal O_{\ul x}$. Multiplying $a_k$ by an even power of $\pi_i$, we can achieve that $v_i(a_k) \geqslant 0$ while preserving the form in a neighborhood of $\ul x$. However, $q$ may be degenerate over $\mathcal O_{\ul\theta}$.

We choose a regular affine open neighborhood $X \subset \widebar X$ of $\ul x$ and $\ul\theta$ such that the coefficients $a_k$ of the form $q$ and the uniformizers $\pi_i$ are regular on $X$ and invertible on $X_\eta$, the reduced closed fiber $X_{v,\red}$ is smooth over $v$, and $\widebar X$ is a compactification of $X$. All these properties except the last can be ensured by choosing $X$ sufficiently small. To ensure the last property as well, we use \cite[Lemma 09NV]{stacks-project}. 

Note that, by construction, each $\pi_i$ is a local parameter at every point of its zero locus in $X$ and $X$ is of constant relative dimension $n \geqslant 1$ over $V$.

\begin{lemma}\label{generic-nondegeneracy}
There exist a finite cover $\rho\colon \mathbf{X} \to X_U$, a quadratic form $\mathbf q$ on $\mathbf{X}$, and a closed subset $Z \subset X$ satisfying the following properties:
    \begin{enumerate}
        \item $\mathbf q$ is nondegenerate and diagonal on $\mathbf{X}$.
        \item $\mathbf q$ is isotropic on $\mathbf{X} - \mathbf{Z}$, where $\mathbf{Z}=\rho^{-1}(Z_U)$
        \item There is a lift $\mathbf\Delta\colon U \to \mathbf X$ of the diagonal $\Delta\colon U \to X_U$ such that $\mathbf q|_{\mathbf\Delta(U)} = q|_U$
        \item The subset $Z$ does not contain any irreducible component of $\widebar X_v$
        \item The nilradical $\nil(\mathbf X_x)$ of each closed fiber $\mathbf X_x$ is generically principal
        \item $\mathbf X$ is irreducible, regular, and admits a compactification
        \item $\mathbf X$ is smooth over $U$ along $\mathbf\Delta(U)$.
    \end{enumerate}
\end{lemma}
\begin{proof}
    It is convenient to consider an auxiliary finite cover $X' \to X$ obtained by adjoining square roots $\sqrt{\pi_i}$ of the $\pi_i$. We define a quadratic form $q'$ on $X'$ as follows:
    \[
        a'_k = a_k \cdot \prod \pi_i^{-v_i(a_k)}
        \quad\text{and}\quad
        q' = \sum a_k'T_k^2.
    \]
    Note that the coefficients of $q'$ are defined on $X$ and are invertible, implying nondegeneracy. Indeed, since $X$ is regular, it suffices to check that the $a'_k$ are invertible at all height-one primes. Since $v_i(a'_k) = 0$, the element $a'_k$ is invertible at $\theta_i$, while invertibility at the remaining primes follows from the invertibility of $a_k$ and $\pi_i$ on $X_\eta$.

    By assumption, $q$ is rationally isotropic, i.e., there exists a rational isotropic unimodular vector $t = (t_1, \dots, t_r)$. Clearly, $t' = (\dots, t_k \cdot \prod (\sqrt{\pi_i})^{v_i(a_k)}, \dots)$ is a $q'$-isotropic vector over $\rat(X')$. Note that the preimage of $\theta_i$ in $X'$ may consist of several points, and $\sqrt{\pi_i}$ is a uniformizer of the local ring at each of these points. Moreover, the valuations of $t_k$ at these points coincide. Multiplying $t'$ by suitable powers $(\sqrt{\pi_i})^{\varepsilon_i}$, where $\varepsilon_i \in \mathbb Z$, we obtain a $q'$-isotropic unimodular vector over $\mathcal O_{X', \ul\theta}$.

    Let $Z' \subset X'$ be the closed subset of points where the vector $t' \prod \sqrt{\pi_i}^{\varepsilon_i}$ fails to be regular or unimodular, and let $Z$ be the image of $Z'$ in $X$. Thus, $t' \prod \sqrt{\pi_i}^{\varepsilon_i}$ is a unimodular vector on $X'-Z'$. By construction, $Z$ does not contain any irreducible component of $\widebar X_v$, which proves (4).

    We now construct $\mathbf{X}$ by adjoining square roots $\lambda_i$ of the regular functions $\pi_i \otimes \pi_i^{-1}$ on $X_U = X \times_V U$ for each $i$ such that $\ul x \cap \widebar{\{\theta_i\}} = \varnothing$. We define a quadratic form $\mathbf q$ analogously to $q'$:
    \[
        \mathbf a_k = (a_k\otimes 1) \cdot \prod \lambda_i^{-2v_i(a_k)}
        \quad\text{and}\quad
        \mathbf q = \sum \mathbf a_kT_k^2,
    \]
    and similarly define the $\mathbf q$-isotropic vector $\mathbf t = (\dots, t_k \cdot \prod \lambda_i^{v_i(a_k)}, \dots)$. We claim that the vector $\mathbf t \prod \lambda_i^{\varepsilon_i}$ is unimodular on $\mathbf X - \mathbf Z$. Both the nondegeneracy of $\mathbf q$ and the unimodularity of this vector can be verified after a faithfully flat base change.

    Let $U' = U \times_X X'$. There is a finite faithfully flat $X_U$-morphism
    \[
        X'_{U'} \to \mathbf X,
        \qquad
        \lambda_i \mapsto \sqrt{\pi_i} \otimes \sqrt{\pi_i}^{-1}.
    \]
    The pullback of $\mathbf q$ is isometric to the form $q' \boxtimes \langle1\rangle$ on $X'_{U'}$, proving (1) (of course, nondegeneracy could also be proved directly). Moreover, the pullback of the vector $\mathbf t \prod \lambda_i^{\varepsilon_i}$ is associated with the vector $t' \prod \sqrt{\pi_i}^{\varepsilon_i}$, proving (2).

    Since $\left.(\pi_i \otimes \pi_i^{-1})\right|_{\Delta(U)} = 1$, there is a lift of $\Delta$ defined by $\lambda_i = 1$. Property (3) is then immediate, since the coefficients of $\mathbf q$ differ from those of $q$ by powers of $\lambda_i$.

    Let us verify property (5). First, the nilradical of $X_v$ is generated by $\pi_i$ at the point $\theta_i$. Since etale base change preserves reducedness, the scheme $X_{v,\red} \times_v x$ is reduced. By Lemma \ref{fiber-product-nilradical}, the nilradical of $X_x$, $x \in \ul x$, is also generated by $\pi_i$ at the points lying over $\theta_i$. We claim that the nilradical of $\mathbf X_x$ is generated by $\lambda_i$ at the points lying over $\theta_i$ whenever $\ul x \cap \widebar{\{\theta_i\}} = \varnothing$, and by $\pi_i$ at the points lying over the remaining $\theta_i$. This follows again from Lemma \ref{fiber-product-nilradical} and the fact that adjoining a root of a generator of the nilradical makes that root a generator of the new nilradical.

    Let us verify property (6). Clearly, $\mathbf X$ is irreducible. Note that $X_U$ is regular, being obtained from the regular scheme $X$ by an essentially smooth base change. To prove the regularity of $\mathbf X$, it suffices to check that $\pi_i \otimes \pi_i^{-1}$ is a local parameter at each point of its zero locus in $X_U$. This property is preserved under essentially smooth base change, $\pi_i$ has this property on $X$, and $\pi_i \otimes \pi_i^{-1}$ is associated with $\pi_i \otimes 1$. To construct a compactification of $\mathbf X$, we use Lemma \ref{adm-comp-vs-normal}. Indeed, $X_U$ admits a compactification $\widebar X_U$, and $\mathbf X$ is finite and surjective over $X_U$.

    Finally, since $X$ is smooth over $V$ at the points of $\ul x$, the scheme $X_U$ is smooth over $U$ along $\Delta(U)$. Since $\pi_i \otimes \pi_i^{-1}$ is invertible along $\Delta(U)$, the scheme $\mathbf X$ is smooth over $U$ along $\mathbf\Delta(U)$. This proves (7).
\end{proof}

\begin{corollary}\label{equating}
    In Lemma \ref{generic-nondegeneracy}, one can additionally
    arrange that $\mathbf q$ is isometric to
    $\langle 1 \rangle \boxtimes q$.
\end{corollary}
\begin{proof}
    Consider the finite etale cover of $\mathbf X$ obtained by adjoining square roots of the units $\mathbf a_k(1 \otimes a_k^{-1})$. By construction, the forms $\mathbf q$ and $\langle 1 \rangle \boxtimes q$ become isometric over this cover.

    Since $\mathbf a_k(1 \otimes a_k^{-1})|_{\mathbf\Delta(U)}=1$, the diagonal $\mathbf\Delta$ admits a lift to this cover, so property (3) is preserved. Also this cover is irreducible. All the remaining properties of Lemma \ref{generic-nondegeneracy} are preserved under any finite etale base change.
\end{proof}


We apply geometric Lemma \ref{geom-lemma} to the $U$-scheme $\mathbf X$ constructed in Corollary \ref{equating}. The lemma yields closed subschemes $D_0$ and $D_1$ of $\mathbf X$ that are finite, flat and one-generated over $U$. Moreover, the form $\langle 1 \rangle \boxtimes q$ is isotropic over $\mathbf X - \mathbf Z$ and, in particular, over $D_0$ and $D_1$. The degrees of $D_0$ and $D_1$ are differ by exactly 1. Hence, one of them is of the constant odd degree. Then by Theorem \ref{springer}, the form $q$ is isotropic.

\bibliographystyle{alpha}
\bibliography{ref}

\end{document}